# Another approach to decide on real root existence for univariate Polynomials, and a multivariate extension for 3-SAT

Deepak Ponvel Chermakani, IEEE Member
deepakc@pmail.ntu.edu.sg   deepakc@myfastmail.com   deepakc@usa.com   deepak.chermakani@ust-global.com

*Abstract:* -    We present six Theorems on the univariate real Polynomial, using which we develop a new algorithm for deciding the existence of atleast one real root for univariate integer Polynomials. Our algorithm outputs that no positive real root exists, if and only if, the given Polynomial is a factor of a real Polynomial with positive coefficients. Next, we define a transformation that transforms any instance of 3-SAT into a multivariate real Polynomial with positive coefficients, if and only if, the instance is not satisfiable.

## 1.   Introduction

One can decide within polynomial time, whether or not a univariate Polynomial with integer coefficients has a real root [1]. However, as the number of variables in the Polynomial increase, the complexity of deciding on existence of a real root increases exponentially with the number of variables. In this paper, we shall develop a new approach for deciding on existence of a real root in the univariate case, and shall analyze its extension to the multivariate case of 3-SAT. We denote an "Integer Polynomial" as a Polynomial having integer coefficients, and a "Real Polynomial" as a Polynomial having real coefficients.

   In the next Section 2, we shall state six theorems on the univariate real Polynomial, using which we shall then construct an algorithm in Section 3, for deciding on existence of a real root in a univariate integer Polynomial. In Section 4, we shall define a transformation that transforms any instance of 3-SAT into a multivariate real Polynomial with positive coefficients, if and only if, the instance is not satisfiable.

## 2.   Six Theorems on the univariate real Polynomial

The given univariate integer Polynomial, for which we are trying to decide on existence of a real root, is of the following form: $P(X) = a_0 + a_1X + a_2X^2 + a_3X^3 + \ldots a_NX^N$, where N is P's degree, where the set $\{a_0, a_1, \ldots a_N\}$ belongs to the set of Integers, where M is the maximum magnitude of elements in $\{a_0, a_1, \ldots a_N\}$, and finally where β represents either the lower bound of the magnitude of the non-zero imaginary component, or the lower bound of the magnitude of the non-zero ratio between the imaginary component and corresponding real component, in the complex roots of P(X).

**Theorem-1**: There exists a real Polynomial R(X) of degree equal to (1 + INT (π p/q)), with positive coefficients, which when multiplied with the Quadratic $((X - p)^2 + q^2)$ where p and q are positive real numbers, gives a resultant real Polynomial having positive coefficients.
**Proof**: We are trying to find the required degree D, of the polynomial $R(X) = X^D + C_1X^{D-1} + C_2X^{D-2} + C_3X^{D-3} + \ldots C_D$, with real coefficients, which when multiplied with $((X - p)^2 + q^2)$, gives a resultant Polynomial with positive coefficients. We proceed by proving 2 Lemmas.
<u>Lemma-1.1</u>: This required degree D, of R(X) is independent of the individual numerical values of p and q, and only depends on the ratio p:q. The higher this ratio, the greater is D
<u>Proof</u>: Consider the product of the quadratic $(X^2 - 2pX + (p^2 + q^2))$ with R (X). If $<1, r_1, r_2 \ldots r_{D+2}>$ depicts the coefficients of $<X^{D+2}, X^{D+1}, X^D, \ldots, X^3, X^2, X^1, 1>$, in this product, then we have:
$r_1 = C_1 - 2p$
$r_2 = C_2 - 2pC_1 + (p^2 + q^2)$
$r_3 = C_3 - 2pC_2 + C_1(p^2 + q^2)$
…
$r_i = C_i - 2pC_{i-1} + C_{i-2}(p^2 + q^2)$
…
$r_D = C_D - 2pC_{D-1} + C_{D-2}(p^2 + q^2)$
$r_{D+1} = -2pC_D + C_{D-1}(p^2 + q^2)$
$r_{D+2} = C_D(p^2 + q^2)$
   Since each of $\{1, r_1, \ldots r_{D+2}\}$ needs to be > 0 by a small amount, then as $r_i \to 0$, for all i in [1,N], we can obtain the following (by starting from $C_1$ and iteratively substituting the value of $C_{i-2}$ and $C_{i-1}$ in the subsequent value of $C_i$):
$C_1 = 2p + \mu_1$
$C_2 = (2p)^2 - (p^2 + q^2) + \mu_2$

$C_3 = (2p)^3 - 2(2p)(p^2 + q^2) + \mu_3$
$C_4 = (2p)^4 - 3(2p)^2 (p^2 + q^2) + (p^2 + q^2)^2 + \mu_4$
$C_5 = (2p)^5 - 4(2p)^3 (p^2 + q^2) + 3(2p)(p^2 + q^2)^2 + \mu_5$

In the previously mentioned sequence of $C_i$, the quantities $\mu_1, \mu_2, \mu_3$, etc, are positive and → 0.

If we denote the ratio of p:q as h, then the previously mentioned sequence of $C_i$ becomes:

$C_1 = 2p + \Delta_1$
$C_2 = p^2 ((2)^2 - (1^2 + h^2)) + \Delta_2$
$C_3 = p^3 ((2)^3 - 2(2)(1^2 + h^2)) + \Delta_3$
$C_4 = p^4 ((2)^4 - 3(2)^2 (1^2 + h^2) + (1^2 + h^2)^2) + \Delta_4$
$C_5 = p^5 ((2)^5 - 4(2)^3 (1^2 + h^2) + 3(2)(1^2 + h^2)^2) + \Delta_5$
…

where the quantities $\Delta_1, \Delta_2, \Delta_3, \Delta_4$, etc, are positive and → 0.

By induction, it is straightforward to obtain the result that the sign of $C_i$ is independent of the individual values of p or q, but only depends on the ratio p:q. And what is actually required is that $C_D$ is just negative, which will mean that D is the degree of R(X) that is just sufficient to ensure that all the coefficients of the resultant polynomial are just sufficiently above zero. The reason is that if $C_D$ is just negative, then we can force $C_D = 0$, making $r_{D+2} = 0$, and $r_{D+1}$ = some positive value because $C_{D-1}$ is still positive, and finally $r_{D-2}$ is obviously some small positive quantity. Thus, D depends only on the ratio p:q. Further, it is obvious that suppose if we were to force q=0, then there can exist no expression R(X) with a finite degree, that can force the product of R(X) and $(X^2 - 2pX + p^2)$ to have zero sign changes in its coefficients, because by Descartes Rule of Signs, we know that $(X^2 - 2pX + p^2)$ has 2 real roots, and therefore the product of R(X) with $(X^2 - 2pX + p^2)$ must also have atleast 2 sign changes while reading its coefficients. Then, as the ratio p:q keeps increasing; the required value of D also increases. Hence Proved Lemma-1.1

Lemma-1.2: The required degree D of R(X) is bounded by the product of the ratio p:q and the irrational π (which is 3.141592…).
Proof: To prove this, we first look at the sequence $S_i$ obtained while deriving the value of π [2]. Next, we will show that for the case of a=p=½, the terms of our sequence $C_i$ (that we introduced in Lemma-1.1) decreases faster than $S_i$, and that $C_1 = S_1$, and also that $C_2 < S_2$. Finally, we invoke Lemma-1.1, to show that sequence $C_i$ decreases faster than $S_i$ irrespective of the value of p.

Consider a circle of radius a, centered at (0, a), where we are trying to inscribe triangles in the right half of the circle wrt Y-axis. Let us consider 2 such triangles as shown in the below figure-1.

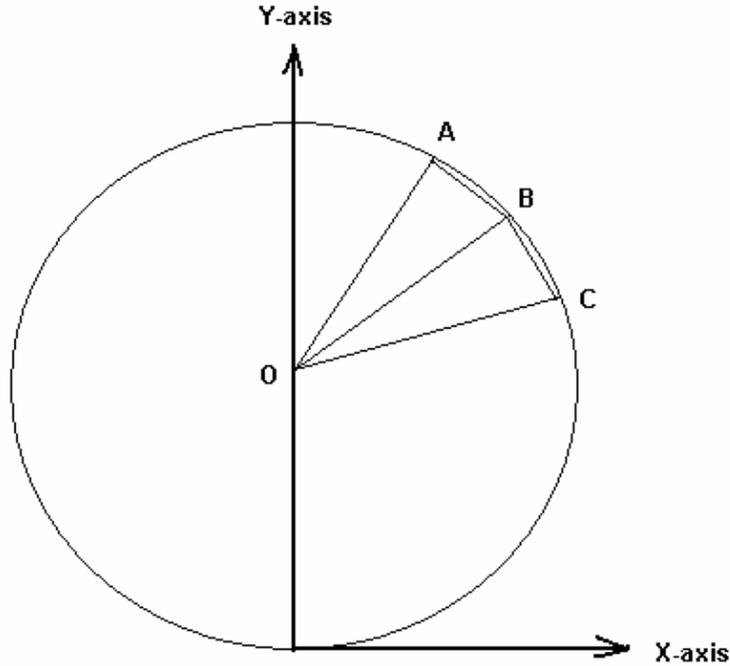

Fig 1.   Two inscribed triangles being considered

In the above figure-1, angles AOB and BOC are denoted to be equal to $2\partial$, and the angle between OA and Y-axis is denoted to be equal to ℮. Let us denote $S_i$, $S_{i+1}$, and $S_{i+2}$ as the lengths of the projections of OA, OB, and OC respectively on the Y-axis, and let us also denote the lengths of segments AB=BC=b/2. Then we obtain 4 equations:

$\sin (\partial) = b / 4a$
$S_i - a = a \cos (℮)$
$S_{i+1} - a = a \cos (℮ + 2\partial)$
$S_{i+2} - a = a \cos (℮ + 4\partial)$

Simplifying the previous 4 equations yields the relationship: $S_{i+2} = 2 S_{i+1} (1 – b^2/(8a^2)) – S_i + b^2/(4a)$

If we put $a=½$, and simplify then it becomes $S_{i+2} = 2 S_{i+1} (1 – b^2/2) – S_i + b^2/2$

Let us revisit the derivation of the irrational $\pi$. As we go on inscribing more and more small triangles, where the top most triangle touches X=1, and the bottom most triangle touches X=0, then as b→0, the number of such triangles being 'n', then using the formula for the length of half of the circumference of a circle, we can say that $nb/2 = \pi/2$, or $n=\pi/b$. We have $S_1=1$ (since the Y-projection of the tip of the top most triangle meets the Y-axis at 1).

Looking back at our original sequence, since $C_i – 2pC_{i-1} + C_{i-2} (p^2 + q^2)$ is almost zero, then putting $p=½$ and replacing i with i+2, we get $C_{i+2} = C_{i+1} – C_i(q^2 + ¼)$, and we also get $C_1 = 1$ (which is equal to $S_1$), and $C_2 = ¾ – q^2$ (which is lesser than $S_2$, because $S_2$ is supposed to be very close to $S_1=1$, when b→0).

It is easy to see that when $C_i$ and $S_i$ are both positive, then $C_i > C_{i+1} > C_{i+2}$, and $S_i > S_{i+1} > S_{i+2}$. Now, we will go on to prove that the rate at which $C_i$ decreases is faster than the rate at which $S_i$ decreases, which, coupled with the fact that $S_i$ reaches zero when $i > \pi/b$, would enable us to conclude that $C_i$ will also become negative before $i > \pi/q$, as q→0. This proof is as follows.

We have the 2 sequences, when $a=p=½$:

$S_{i+2} = 2 S_{i+1} (1 – b^2/2) – S_i + b^2/2$, being the first sequence, and,

$C_{i+2} = C_{i+1} – C_i(b^2 + ¼)$, being the second sequence. (note that we are replacing q with b, as both quantities → 0)

Adding and subtracting terms, we can say that $S_{i+2} = S_{i+1} – S_i(b^2 + ¼) + (1-b^2)S_{i+1} – S_i(¾ - b^2) + b^2/2$. Now $S_i$ decreases at a slower rate than $C_i$, if $(1-b^2)S_{i+1} – S_i(¾ - b^2) + b^2/2 > 0$, which is true if $S_{i+1} – ¾ S_i + b^2(½ – S_{i+1} + S_i) > 0$. It is easy to see that it is true, because the coefficient of $b^2$ is positive, and because the value of $S_{i+1} – ¾ S_i$ is also positive, as $S_{i+1}$ is lesser than $S_i$ by a quantity less than b which we have assumed to be very small. Thus, $C_i$ becomes negative before $i > INT(\pi/b)$, or $INT(\pi/q)$.

In our entire argument above, we have assumed that $p = ½$. However, even if $p \neq ½$, we may recall from Lemma-1.1, that the event of $C_i$ becoming just negative does not depend on p or q, but only on the ratio p:q, hence we can directly conclude that $C_i$ will become negative before $i > \pi (p/q)$. <u>Hence Proved Lemma-1.2</u>

Lemma-1.1 and Lemma-1.2, complete the proof of Theorem-1. However, for verification, we have pasted computer simulations in Table-1, regarding when sequence $C_i$ becomes negative, and it is interesting to see that $(i/(p/q))$ does tend to $\pi$.

Finally, it is also obvious that not only does the resulting Polynomial (product of the Quadratic with R(X)) have positive coefficients, but also R(X) itself has positive coefficients.

**Hence Proved Theorem-1**

Table 1. Computer simulation results for checking when $C_i$ just becomes negative

| p/q | i at which $C_i$ just becomes negative |
|---|---|
| 1 | 3 |
| 2 | 6 |
| 3 | 9 |
| 4 | 12 |
| 6 | 19 |
| 8 | 25 |
| 9 | 28 |
| 12 | 37 |
| 13 | 40 |
| 15 | 47 |
| 19 | 59 |
| 20 | 62 |
| 30 | 94 |
| 34 | 106 |
| 36 | 113 |
| 40 | 125 |
| 41 | 128 |
| 49 | 153 |
| 50 | 157 |
| 51 | 160 |
| 100 | 314 |
| 200 | 628 |
| 1,000 | 3141 |
| 10,000 | 31415 |
| 100,000 | 314159 |
| 1,000,000 | 3141592 |

**Theorem-2**: **Let Q(X) be a real polynomial. Deciding on real root existence for the given univariate polynomial Q(X), is equivalent to deciding the same for $((X – p_1)^2 + q_1^2) ((X – p_2)^2 + q_2^2) ((X – p_3)^2 + q_3^2)…((X – p_N)^2 + q_N^2) = 0$, where $\{p_1, p_2, … p_N\}$ belongs to the set of real numbers, and where $\{q_1, q_2, … q_N\}$ belongs to the set of non-negative real numbers.**

**Proof**: Whether or not Q(X)=0 is solvable in a real interval, is equivalent to whether or not $Q^2(X)=0$ is solvable in that same real interval. Next, we know that the real roots of a polynomial may or may not be repeated, and that the complex roots must occur in conjugate pairs. Therefore, when we square Q(X), we obtain a product of N quadratic components: $((X – p_1)^2 + q_1^2) ((X – p_2)^2 + q_2^2) ((X – p_3)^2 + q_3^2)….. ((X – p_N)^2 + q_N^2)$, where $\{p_1, p_2, … p_N\}$ belongs to the set of real numbers, and where $\{q_1, q_2, … q_N\}$ belongs to the set of non-negative real numbers.
**Hence Proved Theorem-2**

**Theorem-3**: Let Q(X) be a real polynomial. In those quadratic components of $Q^2(X)$, $((X – p_i)^2 + q_i^2)$ where $p_i$ and $q_i$ are positive real numbers, the magnitude of the ratio $p_i:q_i$ is bounded below $(2 \beta^{-1} (MN)^2)$.
**Proof**: Let us denote one of the factor polynomials of $P^2(X)$, to be $P_1(X)$, a polynomial with real coefficients. We denote K as the degree of $P^2(X)$, and L as the maximum coefficient of $P^2(X)$. Then it is clear that K=2N, and L< $NM^2$.

We also define "coefficient normalization" of a polynomial, as the multiplication of every coefficient by a real, such that the smallest non-zero magnitude of coefficients in the polynomial becomes one. Then we state the following 2 lemmas:
<u>Lemma-3.1</u>: After coefficient normalization of $P_1(X)$, the magnitude of every coefficient in $P_1(X)$, is bounded below (KL).
<u>Proof</u>: In a "meta" level, the above Lemma means: - that if $P_1(X)$ has a very high level of "expression precision", then $P_1(X)$ cannot be a factor of $P^2(X)$ that has a very low level of "expression precision", unless the degree of $P^2(X)$ is very high. Example: a high "expression precision" quadratic like $(10101+ 40332X + 809\sqrt{3}X^2)$ cannot be a factor of a low "expression precision" polynomial like $(1 + 4X - 8X^2 + 3X^3 - 2\sqrt{3}X^4+ 3X^5+ 7X^6)$. However, our Lemmas indicate that there is a possibility that $(10101 + 40332X + 809\sqrt{3}X^2)$ might be a factor of $(1 + 4X - 8X^2 + 3X^3 - 2X^4 + 5X^5 - 2X^6 + … 3X^{40332})$. The figure 2 below explains.

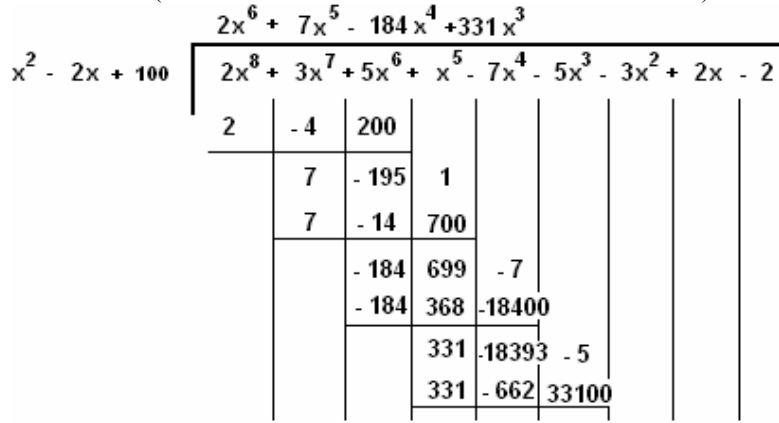

**Fig 2. Example of the uncontrolled increase in the numerical expression of numbers, in the attempted division of a low-precision, low-degree polynomial, by a high-precision quadratic**

We will now verify that (KL) is the upper bound for the magnitude of any coefficient in $P_1(X)$, after coefficient normalization. Consider an example, and let us say $P^2(X) = (1 + 1X + 1X^2 + … 1X^N – 1X^{N+1} – 1X^{N+2} –X^{2N})$. Two factors of $P^2(X)$ are (x-1) and $(1 + 2X + 3X^2 + … (N-1)X^N + (N)X^{N+1} + (N-1)X^{N+2} … 3X^{2N-3} + 1X^{2N-2} + 1X^{2N-1})$. Let $P_1(X)$ be the second factor. This example shows an extreme case, in which the maximum magnitude of coefficients in $P_1(X)$, is equal to (N). Here, K, the degree of $P^2(X)$ is 2N, while L, the maximum coefficient of $P^2(X)$ is 1. Thus (KL) bounds the magnitude of every coefficient in $P_1(X)$. If one tries to alter the factor polynomials by changing the degree, or by changing coefficient magnitudes of the factor polynomials, such that their product is still $P^2(X) = (1 + 1X + 1X^2 + … 1X^N – 1X^{N+1} – 1X^{N+2} –X^{2N})$, then one would find that the maximum magnitude of coefficients of $P_1(X)$ decreases. Similarly, if one tries to change the degree/coefficients of $P^2(X)$, and repeats the experiment with other factor polynomials, it is found that the magnitude of every coefficient of the factor polynomial after coefficient normalization, is bounded by the value of (KL) obtained from the new $P^2(X)$. <u>Hence Proved Lemma-3.1</u>
<u>Lemma-3.2</u>: $(2KL)^{-1}$ < non_zero_magnitude($p_i$) < (KL), and, *Min_Real_Root_Sep* < non_zero_magnitude($q_i$) < (KL)
<u>Proof</u>: In the quadratic $X^2 – 2p_i X – (p_i^2 + q_i^2)$, if non_zero_magnitude($p_i$) is outside the bound [(KL), $(2KL)^{-1}$], or if non_zero_magnitude($q_i$) > (KL), then the quadratic would, after coefficient normalization, have the maximum magnitude of one of its coefficients become greater than (KL). When this happens, it is impossible for us to find a real polynomial, which when multiplied with the quadratic, would yield $P^2(X)$, as per the argument of Lemma-3.1. <u>Hence Proved Lemma-3.2</u>

Obtaining a lower bound on the value of $q_i$ is more difficult. We give a simple argument showing that the behavior of the lower bound of the imaginary component (i.e. $q_i$) of the complex root of $P^2(X)$, is similar to the behavior of the minimum separation between real roots: - Consider the product of (X-r) with (X-(r+Δ)), which yields $X^2 –X(2r+Δ) – (r^2+rΔ)$. Subtracting four times the coefficient of 1 from the square of the coefficient of X yields the quantity $Δ^2$. Compare this with our quadratic $X^2 – 2p_i X – (p_i^2 + q_i^2)$, where subtracting one fourth the coefficient of X from the coefficient of 1 yields the quantity $q_i^2$. And literature suggests that the lower bound for the minimum separation between real roots decreases exponentially with N and M [3][4][5].

That is why we introduced β. Plugging in the values, the magnitude of the ratio $p_i:q_i$ is bounded below $(2 \beta^{-1} (MN)^2)$, if β denotes the smallest magnitude of non-zero imaginary components. But if β denotes the smallest magnitude of non-zero ratios between the imaginary component and the corresponding real component, in the complex roots of P(X), then it is obvious that non_zero_magnitude($p_i:q_i$) is bounded below $\beta^{-1}$. In either case, non_zero_magnitude($p_i:q_i$) is bounded below $(2 \beta^{-1} (MN)^2)$.
**Hence Proved Theorem-3**

**Theorem-4:** Let Q(X) be a real polynomial. Q(X) does not have a positive real root, if and only if, Q(X) is a factor of some real polynomial with positive coefficients.

**Proof:** If Q(X) does not have a positive real root, then $Q^2(X)$ does not have a positive real root, in which case, each of the N quadratic components of $Q^2(X)$, $((X - p_i)^2 + q_i^2)$, can be multiplied with some Polynomial with positive coefficients of degree equal to $(1 + INT (\pi p_i / q_i))$, to produce polynomials with positive coefficients. This means that there exists a real polynomial T(X) with positive coefficients, which when multiplied with $Q^2(X)$, gives a real polynomial V(X) with positive coefficients. The degree of T(X) would be equal to $(1 + SUMMATION (INT (\pi p_i / q_i)))$, over all i as integers in [1,N]. This implies that Q(X) multiplied by Q(X) T(X) would yield V(X). So Q(X) is a factor of V(X). But if Q(X) has a positive real root, then Q(X) cannot be a factor of any real polynomial with positive coefficients, because of Descartes Rule of Signs. Finally, Q(X) cannot be the factor of some real polynomial with positive coefficients, only if Q(X) has a positive real root. Also, Q(X) can be the factor of some real Polynomial with positive coefficients, only if Q(X) does not have a positive real root.
**Hence Proved Theorem-4**

**Theorem-5:** Let Q(X) be a real polynomial. Q(X) does not have a positive real root, if and only if, Q(X) can be multiplied with a real Polynomial U(X) with positive coefficients, to give a real Polynomial with positive coefficients.

**Proof:** If Q(X) does not have a positive real root, then Q(X) may or may not have negative roots, and also may or may not have complex roots with positive real components. We are interested only in those complex roots with positive real components, because these are the ones that cause negative signs to creep into Q(X). As these complex roots occur as conjugates, therefore they may be represented as the product: $((X - p_i)^2 + q_i^2)\ldots((X - p_j)^2 + q_j^2)$. Now it is obvious from Theorem-1, that Q(X) can be multiplied with some U(X) with positive coefficients, to yield a resultant Polynomial with positive coefficients.
**Hence Proved Theorem-5**

**Theorem-6:** Let Q(X) be a real polynomial. If Q(X) has 'k' positive real roots, then Q(X) can be multiplied with a real Polynomial U(X) with positive coefficients to give a real Polynomial with 'k' sign changes in its coefficients.

**Proof:** Q(X) can consist of five types of simplest real factors, the first of the form (X), the second of the form $(X - p_i)$, the third of the form $(X + p_i)$, the fourth of the form $((X - p_i)^2 + q_i^2)$, and the fifth of the form $((X + p_i)^2 + q_i^2)$, where $p_i$ and $q_i$ are positive real numbers. We are interested only in the factors of the second and the fourth types, because these factors are the ones which cause negative signs to creep into Q(X). However, we know from Theorem-5 that the fourth type of factor can be multiplied with a Polynomial to produce a real Polynomial with positive coefficients. This means that all types of factors, except factors of the second type, can be multiplied by some real Polynomial with positive coefficients, so as to yield a real Polynomial with positive coefficients (let us call this Polynomial as B). Say there are 'k' such factors of the second type, i.e. $(X - a_1)(X - a_2)\ldots(X - a_k)$. Let us say that B has 'b' monomials. Now it is obvious that the first factor $(X - a_1)$ can be multiplied with a real Polynomial with positive coefficients so as to yield the Polynomial $(X^b - a_1^b)$. If we multiply B with $(X^b - a_1^b)$, we get another Polynomial $B_1$ having 'b' +ive monomial terms followed by 'b' –ive monomial terms. In other words, $B_1$ has 1 sign change. Again, the second factor $(X - a_2)$ can be multiplied with a real Polynomial with positive coefficients so as to yield the Polynomial $(X^b - a_2^b)$. If we multiply $B_1$ with $(X^b - a_2^b)$, we get another Polynomial $B_2$ having 'b' +ive monomial terms followed by 'b' –ive monomial terms followed by 'b' +ive monomial terms. In other words, $B_2$ has 2 sign changes. It follows by induction that after multiplying $B_{k-1}$ with $(X^b - a_k^b)$, we get Polynomial $B_k$ having 'k' sign changes.
**Hence Proved Theorem-6**

## 3. The new Algorithm for deciding whether a univariate integer Polynomial has a positive real root

The basic idea of our Algorithm stems from Theorem-5, which is to check whether there exists a Polynomial T(X) of a particular degree, which when multiplied with the given integer Polynomial P(X), of degree N and maximum magnitude of coefficients M, gives a resultant polynomial with positive coefficients. To obtain the degree of T(X), we look at Theorem-3, which said that the maximum value that the magnitude of the ratio $p_i:q_i$ can take, if $p_i > 0$ and $q_i \neq 0$, in any quadratic component of $P^2(X)$, is the value $(2 \beta^{-1}(MN)^2)$, which, according to Theorem-1, would need a Polynomial of degree equal to $(1 + INT(2\pi \beta^{-1} (MN)^2))$. Next, there are N such potential quadratic components in $P^2(X)$, so the degree of T(X) is conveniently chosen as $(1 + INT(2\pi \beta^{-1} M^2 N^3))$, because it is harmless to use T(X) of degree higher than required. So, here is our algorithm.
**Start**
**Step-1:** Set V(X) = P(X) T(X), where $T(X) = X^D + T_1 X^{D-1} + T_2 X^{D-2} + \ldots T_D$, and where $D = (1 + INT(2\pi \beta^{-1} M^2 N^3))$
**Step-2:** Use a Linear Programming Solver to check whether or not there exists a set of positive reals $\{T_1, T_2, T_3, \ldots T_D\}$, such that every coefficient of V(X) is greater than zero
**Step-3:** If the answer from Step-2 is YES, then P(X) does not have a real root in $]0, \infty[$, and if the answer from Step-2 is NO, then P(X) has a real root in $]0, \infty[$
**Stop**

In our algorithm, the Linear Programming Solver receives $(1 + INT(2\pi \beta^{-1} M^2 N^3) + N)$ inequations, with integer coefficients whose magnitude is limited to M, and having $(1 + INT(2\pi \beta^{-1} M^2 N^3))$ variables. To check for real root existence in $]0, -\infty[$, we need to simply repeat the algorithm with P(-X). It is trivial to check for real root existence at X=0.

# 4. Application to the multivariate Polynomial derived from 3-SAT

Given an instance of 3-SAT having u binary variables, we can obtain a u-variate integer Polynomial Q, such that Q has a real root in $S_u$, the region defined by positive real values of all the variables $X_1…X_u$, if and only if, the 3-SAT instance is satisfiable. Also, if the 3-SAT instance is not satisfiable, then Q will not have a real root. One such method for obtaining Q is as follows.
Step-1: Take one real variable for each Boolean variable in the given 3-SAT instance.
Step-2: Express each clause by a Polynomial. If a clause involves variables $X_1$, $X_2$, and $X_3$, express it by the Polynomial
$(X_1 - i)^2 (X_2 - j)^2 (X_3 - k)^2$, where i is 1 if $X_1$ appears negated and 2 otherwise, j is 1 if $X_2$ appears negated and 2 otherwise, and k is 1 if $X_3$ appears negated and 2 otherwise.
Step-3: Set Q = SUMMATION (the Polynomial expressing each clause in Step 2, over all clauses in the 3-SAT instance).

We are unable to extend the definition of β to Q, because there can exist complex roots with arbitrarily small imaginary components, if Q has a real root. Even if Q does not have a real root, one can assign a complex number (with an arbitrary small imaginary component) to one of the variables, and find the complex solution set for the other variables, such that Q=0, showing that Q can have a complex root with an arbitrary small imaginary component, irrespective of whether Q has a real root or not.
But it is important to note that if Q does not have a real root, then no point in the real-sketch of Q can come arbitrarily close to the Y=0 Plane. Thus, even though we are unable to define β, we are not stopped from proving the following Theorem for Q.

**Theorem-7**: Let Q be the u-variate Polynomial derived from 3-SAT, and $P_i = (X_i – 1)^2 (X_i – 2)^2$ for all i as integers in [1,u]. There exists a real Polynomial P with positive coefficients, and there exist real Polynomials K, $K_1$, $K_2$ …$K_u$, such that P = Q K + $P_1 K_1$ + $P_2 K_2$ + $P_3 K_3$ + … + $P_u K_u$, if and only if, Q does not have a real root.
**Proof**: We prove this via four Lemmas.
Lemma-7.1: P exists, if Q does not have a real root.
Proof: Let $Q = Q_{N1} X_1^{N1} + Q_{N1-1} X_1^{N1-1} … + Q_2 X_1^2 + Q_1 X_1 + Q_0$, where $Q_{N1} …Q_0$ denote Polynomials in $X_2…X_u$. Actually $N_1$ = 2, but we shall continue using the symbol $N_1$ to make our proof more general. Let R = Q + Q', where Q' = $(X_1^{N1} +1)$ (SUMMATION $((X_i – 1)^2 (X_i – 2)^2$, over all i as integers in [2,u]))$^{MAX\_1}$. Here MAX_1 is the degree of Q. If L and L' denote the families of univariate Polynomials in $X_1$ obtainable from Q and R respectively, by substituting all the other variables as any Real, then one can identify two motivations why we chose Q' in this manner. The first motivation is that we want to ensure that the Question of whether or not R has a real root, is equivalent to the Question of whether or not Q has a real root. In other words, the members of L obtained by substituting $X_i$ as either 1 or 2 for all i as integers in [2,u], should be a subset of the members of L'. The second motivation is that we want to ensure that every member of L', can always be multiplied by a univariate Polynomial $U(X_1)$ of degree say WORST_CASE_MIN, to yield a univariate Polynomial in $X_1$, with positive coefficients. The members of L face a problem in the sense that the degree of $U(X_1)$ can grow unbounded, as evident in the example of the bivariate Polynomial Q = $(X_1 – X_2)^2 + 1$. But the members of L' do not face this problem, and WORST_CASE_MIN does exist (i.e. it is bounded by some natural number) for the members of L', because as the magnitude of one or more of the variables $X_2…X_u$ tends to ∞, then the members of L' take the form $C_{N1} X_1^{N1} + C_{N1-1} X_1^{N1-1} … + C_2 X_1^2 + C_1 X_1 + C_0$, where $C_0$ and $C_N$ are both very large and positive, and their magnitude is much greater than the magnitude of the remaining coefficients $C_{N-1}…C_1$. So $U(X_1)$ can afford to become the Polynomial $(X_1^{N1} + X_1^{N1-1} … + X_1^2 + X_1 + 1)$, in which case it can be seen that the product of this $U(X_1)$ and any member of L' obtained by substituting one or more of $X_2…X_u$ in R as a real whose magnitude is greater than some real number say λ, will have positive real coefficients.
When we consider members of L' that are obtained by substituting $X_2…X_u$ as values whose magnitude is lesser than λ, then we can still show that the $U(X_1)$ can be some Polynomial such that the product of that member of L' with $U(X_1)$ is a Polynomial with positive coefficients. Let us put C=WORST_CASE_MIN, and since WORST_CASE_MIN can be assumed to be greater than $N_1$ without loss of generality, we can rewrite $R = R_0 + R_1 X_1 + R_2 X_1^2 + ... + R_C X_1^C$. From the equation R F = E, we can say that the following (C+D+2) Inequations are being satisfied for all real values of $X_2…X_u$:
$F_0 > 0$,
$F_1 > 0$,
$F_2 > 0$,
…
$F_C > 0$,
$F_0 R_0 > 0$,
$F_1 R_0 + F_0 R_1 > 0$,
$F_2 R_0 + F_1 R_1 + F_0 R_2 > 0$,
$F_3 R_0 + F_2 R_1 + F_1 R_2 + F_0 R_3 > 0$,
…
$F_C R_{C-1} + F_{C-1} R_C > 0$,
$F_C R_C > 0$.
In the above Inequations, each of $F_0…F_C$ are Continuous in $X_2…X_u$, evaluating to positive real values for all real values of $X_2…X_u$. The reason for this is that, if we consider the above (C+D+2) Inequations, we can see that the feasible solution set of the unknowns $F_0…F_C$, are bounded by a dynamic Polytope, the boundaries of which are changing continuously with $R_0…R_C$, that happen to be Polynomial functions of $X_2 \ldots X_u$. It is important to note that though the number of boundaries of this Polytope can vary discretely, still these boundaries are shifting continuously, allowing a continuous solution set of $F_0 \ldots F_C$ to be identified.

It is known that any continuous function within a real interval can be approximated as a Polynomial. So we approximate $F_0 \ldots F_C$ as Polynomials say $P_0 \ldots P_C$ within the region $S = \{-\lambda < X_i < +\lambda$, where $X_i$ is real, for all i as integers in $[2,u]\}$. Using the fact that $F_0 \ldots F_C$ tend to Polynomials within the region S, and the fact that they tend to 1 outside the region S, we can simply say, for all k as integers in $[0,C]$, that $F_k = P_k / (1 + \text{SUMMATION} ((X_i/\lambda)^Z,$ over all i as integers in $[2,u])) + 1 / (1 + \text{SUMMATION} ((\lambda/X_i)^Z,$ over all i as integers in $[2,u]))$, where Z is some large even natural number.

Thus, for all k as integers in $[0,C]$, $F_k$ can be depicted as $num_k/den_k$, where $num_k$ and $den_k$ are Polynomials in $X_2 \ldots X_u$, that evaluate to positive real numbers, for all real values of $X_2 \ldots X_u$. However, if we multiply each of the $(C+D+2)$ Inequations with PRODUCT ($den_k$, over all k as integers in $[0,C]$), then it is possible to identify $F_0 \ldots F_C$ as Polynomials in $X_2 \ldots X_u$, such that $R F = E$, where $R = R_0 + R_1 X_1 + R_2 X_1^2 + \ldots + R_{N1} X_1^{N1}$, $F = F_0 + F_1 X_1 + F_2 X_1^2 + \ldots + F_C X_1^C$, and $E = E_0 + E_1 X_1 + E_2 X_1^2 + \ldots + E_D X_1^D$, where each of $F_0 \ldots F_C$ and each of $E_0 \ldots E_D$ are Polynomials in $X_2 \ldots X_u$, and where each of $F_0 \ldots F_C$ and each of $E_0 \ldots E_D$ evaluate to positive real values for all real values of $X_2 \ldots X_u$.

We have now basically 'eliminated' the variable $X_1$. By 'eliminating $X_1$', we mean that we have expressed the solvability of Q as the solvability of some Polynomial, which when expressed as a function of $X_1$, has all its coefficients positive for all real values of $X_2 \ldots X_u$. Before proceeding to eliminate $X_2$, we need to add $(X_2^{N2} +1)$ (SUMMATION $((X_i - 1)^2 (X_i - 2)^2$, over all i as integers in $[3,u]))^{MAX\_k\_2}$ to $F_k$, for all k as integers in $[0,C]$. Here MAX_k_2 is the degree of the Polynomial $F_k$. Similar to how we eliminated $X_1$, we need to continue, systematically eliminating variables $X_2 \ldots X_{u-1}$. Once that happens, we are finally left with a set of univariate integer Polynomials in $X_u$, after which it is trivial to check using our procedure in Section-3, whether there exist univariate real Polynomials, which when multiplied with each univariate integer Polynomial in the set, will yield univariate Polynomials with positive real coefficients. And, if there exist such Polynomials, then the 3-SAT instance is unsatisfiable, else it is satisfiable.

We shall now focus on proving that Q is unsatisfiable, if and only if, it is possible to identify Polynomial P with positive real coefficients, and real Polynomials K, $K_1$, $K_2$,…$K_u$, such that $P = Q K + P_1 K_1 + P_2 K_2 + P_3 K_3 + \ldots + P_u K_u$, where $P_i = (X_i - 1)^2 (X_i - 2)^2$ for all i as integers in $[1,u]$. For better clarity, let us depict SUM_j_u = SUMMATION $((X_i - 1)^2 (X_i - 2)^2$, over all i as integers in $[j,u])$. The first step of eliminating $X_1$ is the same as described in the previous-to-previous paragraph. In the second step of eliminating $X_2$, it is similar except that instead of adding quantities separately to each $F_k$, we try to get an expression for a quantity that can be directly added to the Polynomial obtained after eliminating $X_1$. And this quantity is $(1+ X_1 + X_1^2 + \ldots + X_1^D) (X_2^{N2} +1)$ SUM_3_u$^{MAX\_2}$, where MAX_2 and $N_2$ are some natural numbers. Continuing on in this way, we can extend our argument mentioned in the above paragraph, to say that Q is unsatisfiable, if and only if, the following u-variate Polynomial has positive real coefficients:

----- Polynomial description starts here -----
$((\ldots u\_minus\_1\_opening\_brackets$
$(Q + (X_1^{N1} +1) \text{SUM\_2\_u}^{MAX\_1})F +$
$(1+ X_1 + X_1^2 + \ldots + X_1^D) (X_2^{N2} +1) \text{SUM\_3\_u}^{MAX\_2})F' +$
$(1+ X_1 + X_1^2 + \ldots + X_1^D) (1+ X_2 + X_2^2 + \ldots + X_2^{D'}) (X_3^{N3} +1) \text{SUM\_4\_u}^{MAX\_3})F'' +$
…
$(1+ X_1 \ldots + X_1^D) (1+ X_2 \ldots + X_2^{D'}) \ldots (1+ X_{u-1} \ldots + X_{u-1}^{D''''u\_minus\_1\_dashes}) (X_{u-1}^{N(u-1)} +1) \text{SUM\_u-1\_u}^{MAX\_u-1})F'''^{u\_minus\_1\_dashes}$
----- Polynomial description ends here -----

In the above description of the Polynomial, $\{F, F', F'' \ldots F'''^{K\_dashes} \ldots F'''^{u\_minus\_1\_dashes}\}$ is a set of some Polynomials, where $F'''^{K\_dashes}$ is a Polynomial in variables $X_{K+1} \ldots X_u$. Also, $\{D,D',D'' \ldots D'''^{u\_minus\_1\_dashes}, N_1, N_2, N_3 \ldots N_{u-1}, MAX\_1, MAX\_2 \ldots MAX\_u-1\}$ is a set of some natural numbers.

If we denote the above Polynomial as P, then P may be re-written as follows. $P = Q K + P_1 K_1 + P_2 K_2 + P_3 K_3 + \ldots + P_u K_u$, where $P_i = (X_i - 1)^2 (X_i - 2)^2$ for all i as integers in $[1,u]$, and where K, $K_1$, $K_2$,…$K_u$ are some real Polynomials.
Hence Proved Lemma-7.1

Lemma-7.2: P does not exist, if Q has a real root.
Proof: If Q has a real root in $S_u$, then it is impossible to identify univariate Polynomials $U(X_1)$ of finite degree, which when multiplied with any member of the L' family, would yield univariate polynomials with positive coefficients.
Hence Proved Lemma-7.2

Lemma-7.3: Q does not have a real root, if P exists
Proof: The evaluation of P is positive at every point in $S_u$, so if P exists, then Q cannot have a real root.
Hence Proved Lemma-7.3

Lemma-7.4: Q does have a real root, if P does not exist
Proof: In our initial definition of Q, we constrained Q to either have a real root in $S_u$, or have no real root in space. We also just proved in Lemma-7.1 that P exists, if Q does not have a real root. Therefore, Q does have a real root in $S_u$, if P does not exist.
Hence Proved Lemma-7.4
**Hence Proved Theorem-7**

## 5.   Two Conjectures that P=NP

Our first Conjecture is with respect to Theorem-5 of our paper, where we can see that a given univariate integer Polynomial P(X) does not have a positive real root, if and only if, there exists a real Polynomial T(X), such that P(X) T(X) is a Polynomial with positive coefficients. As we have argued, the degree of T(X) is likely to increase exponentially with the degree of P(X) because $\beta^{-1}$ is likely to be bounded by an exponential function of the degree of P(X). Next, as described in [6], an instance of 3-SAT can be converted in polynomial time to the question of whether or not some P(X) has a real root. It is also known that determining the solvability of sparse integer Polynomials is NP-hard [7]. In both [6] and [7], the degree of P(X) obtained from 3-SAT, is exponential to the size of the 3-SAT instance. Thus, our Technique described in Section-3, will need time that is doubly exponential to the size of the 3-SAT instance. However, this does not prevent someone from developing a smarter technique to verify the existence of T(X), such that P(X) T(X) will have positive coefficients.

Our second Conjecture is with respect to Theorem-7, where we said that given the u-variate Polynomial Q obtained from 3-SAT, we can obtain a u-variate Polynomial P that will have positive coefficients, if and only if, Q is unsatisfiable. Again, here, the degree and the number of monomial terms in Polynomials K, $K_1$, $K_2$…$K_u$, is likely to increase exponentially with the number of variables in Q. However, this again does not prevent someone, from developing a smarter technique to verify the existence of K, $K_1$, $K_2$…$K_u$, such that P will have positive coefficients.

One inspiration we draw for such smarter techniques, is from the example of SLPs (Straight Line Programs) [6], which give a polynomial time definition for univariate integer Polynomials, though there could be an exponential number of monomial terms in the Polynomial (note that it would take exponential effort to define the Polynomial by enlisting its coefficients).

## 6.   Conclusion

In this paper, we first presented six Theorems on the univariate real Polynomial, using which we presented an algorithm for deciding whether or not an integer Polynomial P(X) has a positive real root. Our Algorithm basically checks whether there exists a real Polynomial T(X) of a particular degree, which when multiplied with P(X), gives a real Polynomial with positive coefficients, in which case the conclusion would be that P(X) has no positive real root, else the conclusion would be that P(X) has a positive real root. The running time of our algorithm is bounded by a polynomial function of N, M and $\beta^{-1}$, where β is smallest magnitude of non-zero imaginary components in the complex roots of P(X).

Next, we found that though the definition of β could not be extended to multivariate Polynomials, we could still present a Theorem which transforms the u-variate integer Polynomial Q derived from 3-SAT, into a u-variate Polynomial P with positive real coefficients, if and only if, the 3-SAT instance is not satisfiable.

Finally, we concluded that a polynomial time approach to 3-SAT can be discovered, if smarter techniques are developed for verifying the existence of 'unknown' univariate or multivariate Polynomials in arithmetic operations involving the 'unknown' and 'known' Polynomials, where the Polynomial that results from the arithmetic operations, has positive real coefficients.

**References**
[1] Basu, Pollack and Roy, *Algorithms in Real Algebraic Geometry*, Second Edition, ISBN 3-540-33098-4
[2] Beckmann, Petr., *A History of PI*, New York: Barnes and Noble Books, 1971
[3] Jaan Kiusalaas, *Numerical Methods in Engineering with MATLAB*, Cambridge University Press, 2005
[4] George E. Collins, *Polynomial Minimum Root Separation*, Journal of Symbolic Computation, 32, pages 467-473, 2001
[5] Siegfried M. Rump, *Mathematics of Computation*, Vol. 33, No. 145, pages 327-336, (Jan., 1979)
[6] Daniel Perrucci, Juan Sabia, *Real root s of univariate Polynomials and Straight Line Programs,* Vol 5, Issue 3, Journal of Discrete Algorithms, pages 471-478, Sep 2007
[7] David A. Plaisted, Theoretical Computer Science 31 (1984), pp. 125-138

**Acknowledgment**
I, Deepak Ponvel Chermakani, am grateful to Dr. David Moews (http://djm.cc/dmoews.html) for the enlightening discussions held with him. I am most grateful to my parents for their sacrifices in bringing me up.

**About the Author**
I, Deepak Ponvel Chermakani, have developed this algorithm and have written this paper out of my own interest and initiative, during my spare time. I am presently a Software Engineer in US Technology Global Private Ltd (www.ust-global.com), where I have been working since Jan 2006. From Feb 2006 to Apr 2006, I attended a 3 month part-time course in *Data Structures and Algorithms* from the Indian Institute of Information Technology and Management Kerala (iiitmk.ac.in), which was sponsored by UST Global. From Jul 1999 to Jul 2003, I studied full-time Bachelor of Engineering degree in *Electrical and Electronic Engineering* from Nanyang Technological University in Singapore (www.ntu.edu.sg). Before that, I completed my full-time high schooling from National Public School in Bangalore in India.

I wish to immediately join a University or Institute as a Masters/PhD student, because I wish to devote myself full-time to study more properties of problems in computational mathematics.